\documentclass[12pt]{amsart}
\usepackage[english]{babel}
\usepackage{amsmath}
\usepackage{amssymb}
\usepackage{amsthm}
\usepackage{yhmath}
\usepackage{url}
\usepackage{graphicx}
\usepackage{caption}
\usepackage{enumitem}
\usepackage{algorithmicx,algpseudocode}
\usepackage[section]{algorithm}
\usepackage{wrapfig}
\usepackage{subcaption}
\usepackage{varwidth}
\usepackage{todonotes}
\newtheorem{theorem}{Theorem}[section]

\theoremstyle{remark}

\newtheorem{example}[theorem]{Example}

\usepackage{csquotes}

\DeclareMathSymbol{\minus} {\mathord}{operators}{"2D}

\usepackage[colorlinks=true, allcolors=blue]{hyperref}

\title{Random graphs from random matrices}

\author{Igor Rivin}
\address{Mathematics Department, Temple University}
\email{rivin@temple.edu}
\date{\today}
\thanks{The author would like to thank the authors of the \texttt{Julia} programming language and of \emph{Wolfram Mathematica}}
\begin{document}
\begin{abstract}
In the paper \cite{giusti2015clique}, the authors  introduced the \emph{order complex} corresponding to a symmetric matrix. In this note, we use it to define a class of models of random graphs, and show some surprising experimental results, showing sharp phase transitions.
\end{abstract}
\maketitle

\section{Introduction}
In the paper \cite{giusti2015clique} the authors introduce the "order complex'' associated to a (symmetric) matrix. Briefly, we view the symmetric $n\times n$ matrix $M$ (with its diagonal set to zero) as the adjacency matrix of the complete graph $K_n,$ and now we produce an increasing family of graphs, starting with the completely disconnected graph on $n$ vertices, and then adding edges in order of increasing size of the corresponding entry of the matrix $M,$ until $p * n(n-1)/2$ edges have been added (in other words, the edge density in the graph is $p$). It is now natural to look at different models of random matrices, use them to generate random graphs, and see what the properties of the random graphs are. 
\begin{example}
\label{ereg}
Suppose $M$ is drawn from the ensemble of symmetric matrices with i.i.d Gaussian entries (note: for this model it is irrelevant what the mean of the Gaussian is). Then the random graphs are nothing but the much studied Erd\"os-R\'enyi random graphs.
\end{example}
\begin{example}
\label{rank1}
In the upcoming paper \cite{curtorivin2019rank} we generate a random vector $v$ and look at the rank one matrix $M(v)=v^t v$ - in the case where the entries of $v$ are iid $\mathcal{N}(0, 1),$ this is a Wishart ensemble. However, if we pick the entries of $v$ to be uniform in $[0, 1],$ we get a model with other properties \footnote{In the paper \cite{curtorivin2019rank} we look at the associated clique complexes, not the graphs per se}.
\end{example}
\begin{example}
\label{pointcloud}
This is, in a way, \emph{the} motivating example: consider a point cloud approximating some shape in $\mathbb{R}^n$ (usually for $n=2, 3$), and let $M$ be the distance matrix of this cloud (that is, the $M_{ij}$ equals the distance between the $i$th and the $j$th points in the cloud.
\end{example}

In this paper we look at the Laplacian eigenvalues of the graphs we construct. There are (at least) two ways to define the Laplacian matrix of a graph.
The first, and simplest is 
\[
L = D - A,
\]
where $D$ is the diagonal matrix of degrees of vertices and $A$ is the adjacency matrix of the graph $G.$
The second is the \emph{normalized symmetric Laplacian}  (see \cite{chung1997spectral}):
\[
\mathcal{L} = D^{-\frac12} L D^{-\frac12},
\]
with $L$ as above.

It turns out that the normalized Laplacian is much better behaved. Note that the normalized Laplacian spectrum is contained between $0$ and $2,$ and the mean is at $1,$ since the trace of $\mathcal{L}$ is always equal to $n.$
\section{Spectral Gap}
\subsection{Erd\"os-R\'enyi model} .
\begin{figure}
  \begin{subfigure}[b]{0.45\textwidth}
    \includegraphics[width=\textwidth]{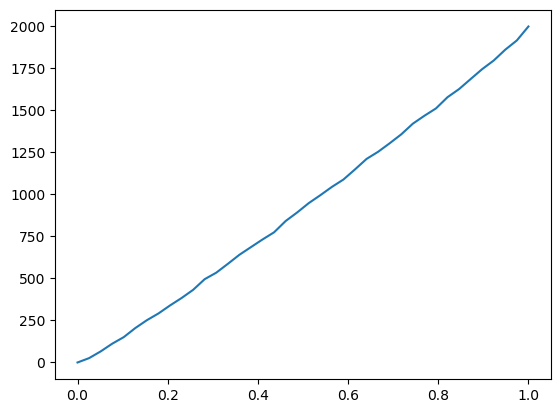}
    \caption{Raw spectral gap}
    \label{fig:0}
  \end{subfigure}
  \begin{subfigure}[b]{0.45\textwidth}
    \includegraphics[width=\textwidth]{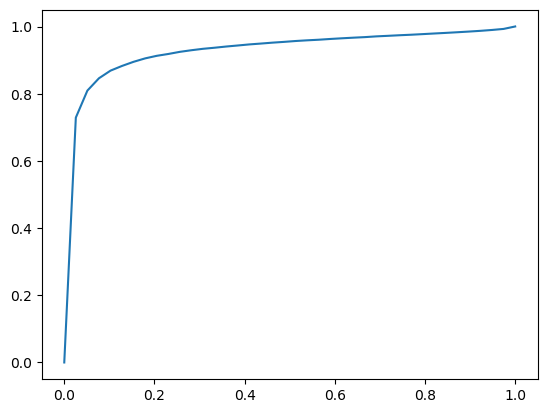}
    \caption{Normalized spectral gaps.}
    \label{fig:1}
  \end{subfigure}
  \caption{ER Model}\label{er}
\end{figure}
We see that the ``raw'' spectral gap - Figure \ref{fig:0} increases linearly from $0$ to the value of the complete graph ($K_2000$ in this case), while the normalized gap - Figure \ref{fig:1} - is asymptotic to $1.$ The latter case \emph{has} been studied - $\lambda_2 \asymp C 1- C n^{-\frac12},$ see \cite{hoffman2012spectral}, but the former seems to be a new observation.
\subsection{Positive rank one model}
\begin{figure}
  \begin{subfigure}[b]{0.45\textwidth}
    \includegraphics[width=\textwidth]{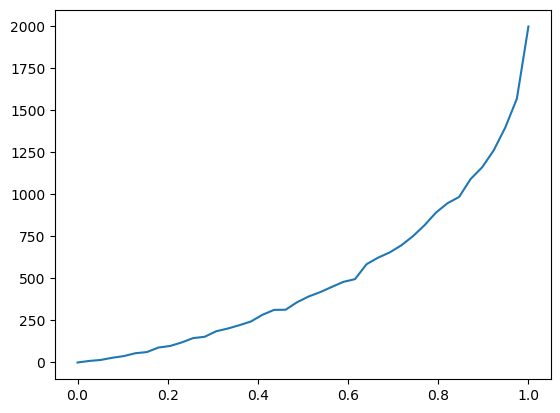}
    \caption{Raw spectral gap}
    \label{fig2:0}
  \end{subfigure}
  \begin{subfigure}[b]{0.45\textwidth}
    \includegraphics[width=\textwidth]{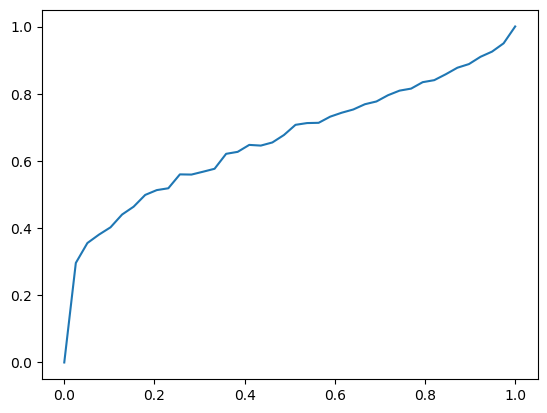}
    \caption{Normalized spectral gap}
    \label{fig2:1}
  \end{subfigure}
  \caption{positive rank 1}\label{posr1}
\end{figure}
We notice that the raw spectral gap - Figure \ref{fig2:0}-  seems to increase like $\sqrt{p},$ while the normalized gap - Figure \ref{fig2:1} - is increasing \emph{linearly} to $1.$ To confirm the first observation, let us plot the square root of the gap:
\begin{figure}
    \centering
    \includegraphics{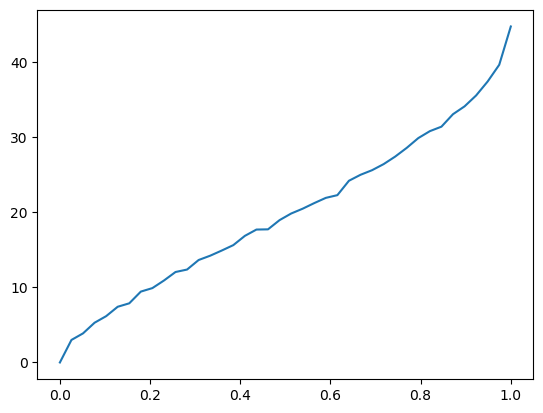}
    \caption{Square root of raw spectral gap for positive rank 1 model}
    \label{fig2:2}
\end{figure}
Note that Figure \ref{fig2:2} is consistent with quadratic growth of the spectral gap. It is also interesting that the two ends (near the completely disconnected and complete graphs) seem symmetric.
\subsection{Rank 1 Wishart model}
\begin{figure}
  \begin{subfigure}[b]{0.45\textwidth}
    \includegraphics[width=\textwidth]{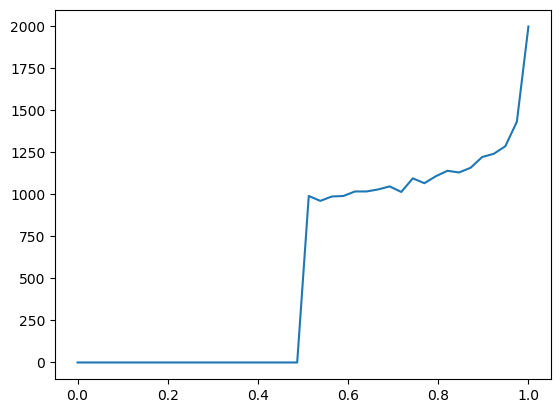}
    \caption{Raw spectral gap}
    \label{fig3:0}
  \end{subfigure}
  \begin{subfigure}[b]{0.45\textwidth}
    \includegraphics[width=\textwidth]{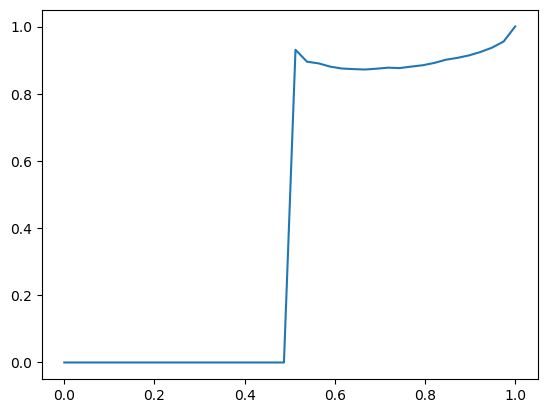}
    \caption{Normalized spectral gaps}
    \label{fig3:1}
  \end{subfigure}
  \caption{Wishart rank 1}\label{wish1}
\end{figure}
The evolution of the spectral gap (see Figure \ref{wish1}) looks starkly different in the Wishart model. Part of the explanation is that (as noted in \cite{curtorivin2019rank}), the graph stays bipartite for low density, until at (roughly) density $p=\frac12$ it becomes complete bipartite (recall that the Laplace eigenvalues of $K_{m, n}$ are $m+n,n, m, 0,$ with multiplicities $1, m-1, n-1, 1.$ However, this explains only some of the features of the evolution (in particular, the sharp phase transition just before the graph becomes complete bipartite and the non-monotonicity of the function).
\subsection{Point clouds}
\begin{figure}
%\label{er}
  \begin{subfigure}[b]{0.45\textwidth}
    \includegraphics[width=\textwidth]{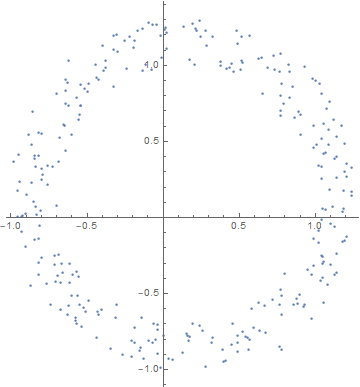}
    \caption{Noisy circle}
    \label{fig4:0}
  \end{subfigure}
  \begin{subfigure}[b]{0.45\textwidth}
    \includegraphics[width=\textwidth]{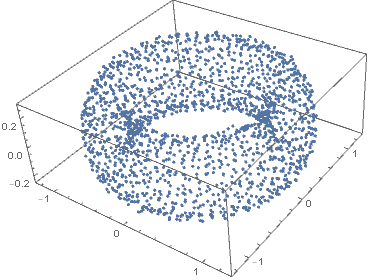}
    \caption{Noisy torus.}
    \label{fig4:1}
  \end{subfigure}
\end{figure}

We now look at the "motivating examples" - point clouds in low-dimensional spaces. The point clouds we look at are the noisy circle and the noisy torus, both found in the \texttt{Eirene} (\cite{henselmanghristl6}) distribution - see Figures \ref{fig4:0} and \ref{fig4:1}.
We convert these point clouds into distance matrices, and see the following spectral behavior:
\begin{figure}
  \begin{subfigure}[b]{0.45\textwidth}
    \includegraphics[width=\textwidth]{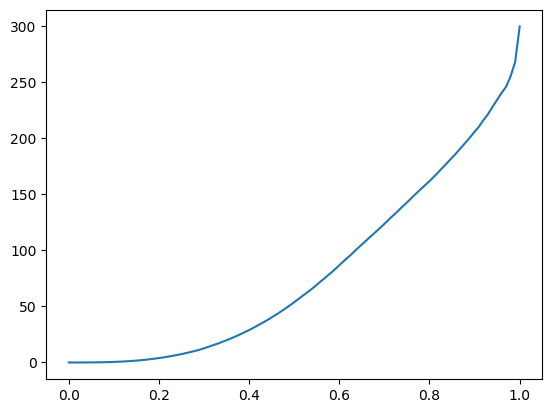}
    \caption{Raw spectral gap}
    \label{fig5:0}
  \end{subfigure}
  \begin{subfigure}[b]{0.45\textwidth}
    \includegraphics[width=\textwidth]{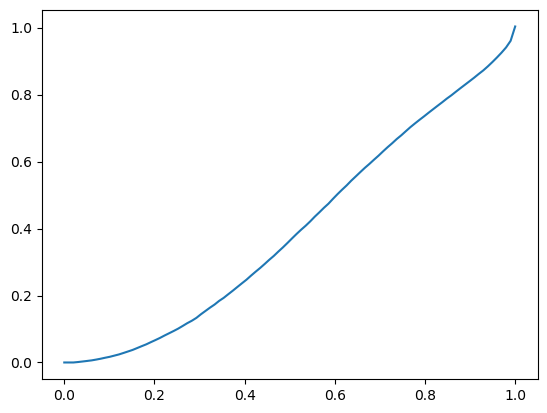}
    \caption{Normalized spectral gap}
    \label{fig5:1}
  \end{subfigure}
  \caption{Noisy circle spectral curves}\label{nc}
\end{figure}
\begin{figure}
  \begin{subfigure}[b]{0.45\textwidth}
    \includegraphics[width=\textwidth]{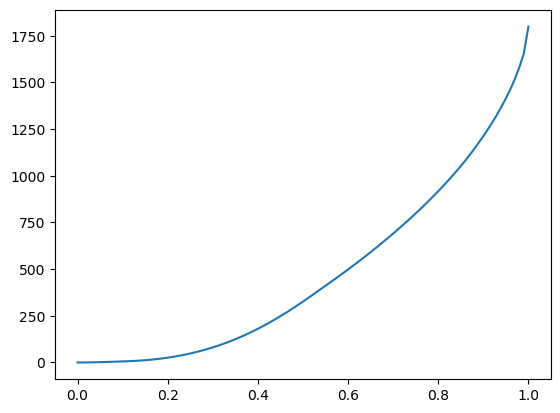}
    \caption{Raw spectral gap}
    \label{fig6:0}
  \end{subfigure}
  \begin{subfigure}[b]{0.45\textwidth}
    \includegraphics[width=\textwidth]{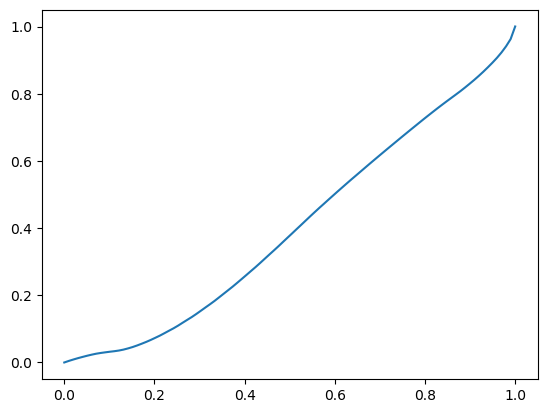}
    \caption{Normalized spectral gap}
    \label{fig6:1}
  \end{subfigure}
  \caption{Noisy torus spectral curves}\label{nt}
\end{figure}
It is quite obvious to the naked eye that the spectral gap curves in Figures \ref{nc} and \ref{nt} are very similar to those in the positive rank one case (Figure \ref{posr1}
\section{Spectral densities}
\subsection{Erd\"os-R\'enyi}
\begin{figure}
  \begin{subfigure}[b]{0.45\textwidth}
    \includegraphics[width=\textwidth]{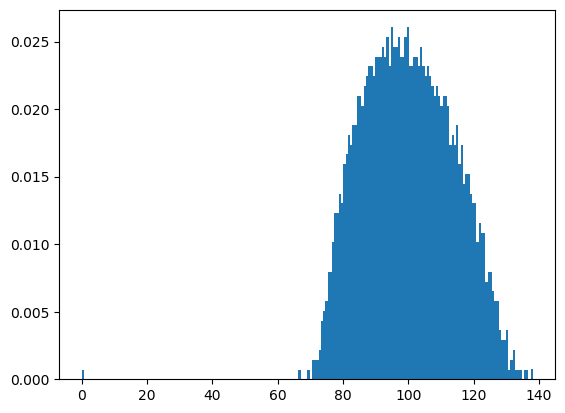}
    \caption{Raw spectral density}
    \label{fig7:0}
  \end{subfigure}
  \begin{subfigure}[b]{0.45\textwidth}
    \includegraphics[width=\textwidth]{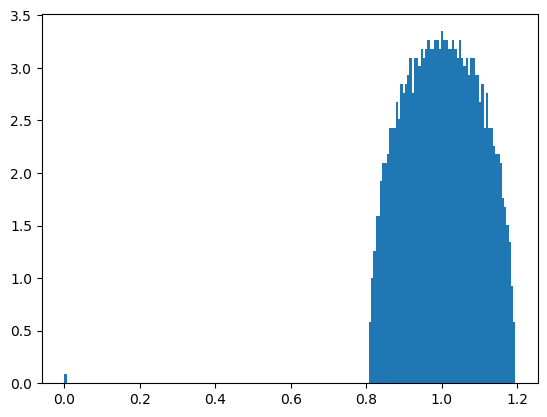}
    \caption{Normalized spectral density}
    \label{fig7:1}
  \end{subfigure}
  \caption{Spectral density at $p=0.05$}\label{goe05}
\end{figure}

\begin{figure}
  \begin{subfigure}[b]{0.45\textwidth}
    \includegraphics[width=\textwidth]{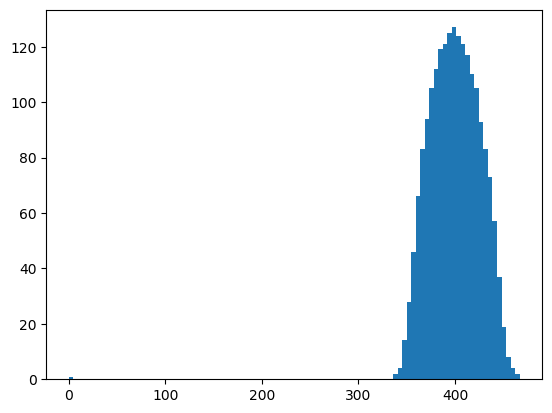}
    \caption{Raw spectral density}
    \label{fig8:0}
  \end{subfigure}
  \begin{subfigure}[b]{0.45\textwidth}
    \includegraphics[width=\textwidth]{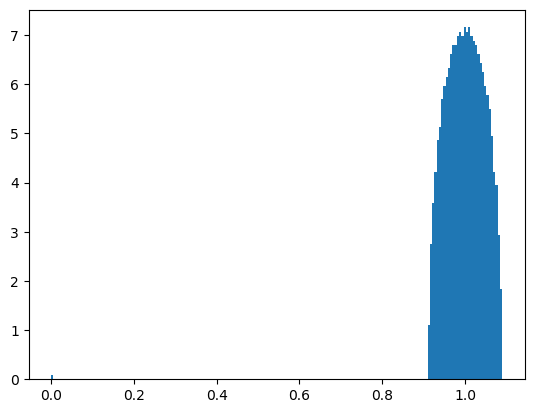}
    \caption{Normalized spectral density}
    \label{fig8:1}
  \end{subfigure}
  \caption{Spectral density at $p=0.2$}\label{goe2}
\end{figure}

The spectral density of the Erd\"os-R\'enyi random graph has been extensively studied (see, for example \cite{erdHos2013spectral}) - the "raw" spectrum seems to have been more extensively studied, and found to satisfy the semicircle law (as the reader might be convinced by looking at the figures \ref{goe05} and \ref{goe2}).

We see that the shapes (whatever that means) of the curves stabilize fairly quickly, and only the width is shrinking with increasing $p.$ it is thus natural to look at the width as a function of $p.$ Instead of the width (which is a little hard to define, we just look at the standard deviation of the empirical distribution of eigenvalues. Let us do it for the normalized spectrum:
\begin{figure}
  \begin{subfigure}[b]{0.45\textwidth}
    \includegraphics[width=\textwidth]{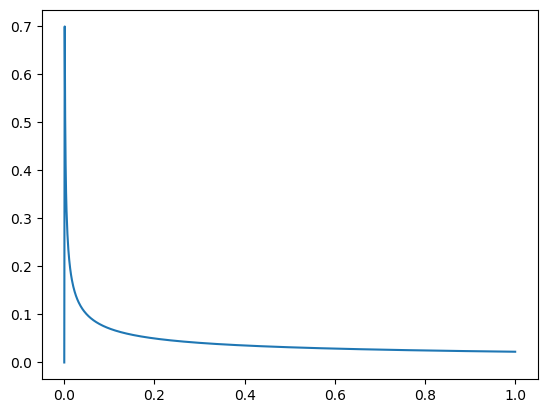}
    \caption{Entire graph}
    \label{fig9:0}
  \end{subfigure}
  \begin{subfigure}[b]{0.45\textwidth}
    \includegraphics[width=\textwidth]{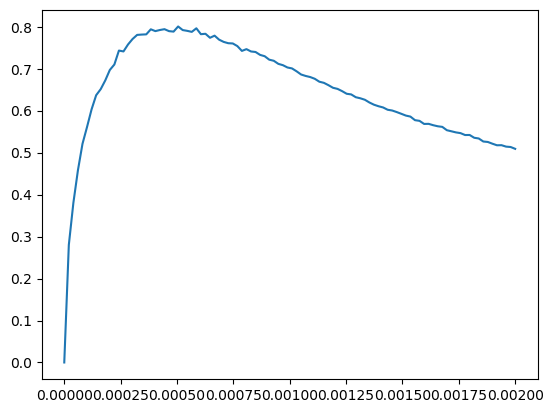}
    \caption{Zoom around maximum}
    \label{fig9:1}
  \end{subfigure}
  \caption{Standard deviation of spectral density}\label{goestd}
\end{figure}
We see in Figure \ref{goestd} that the standard deviation rises sharply until $p=1/n,$ and then declines.
\subsection{Positive rank 1}
First let us look at the spectral distribution:
\begin{figure}
  \begin{subfigure}[b]{0.45\textwidth}
    \includegraphics[width=\textwidth]{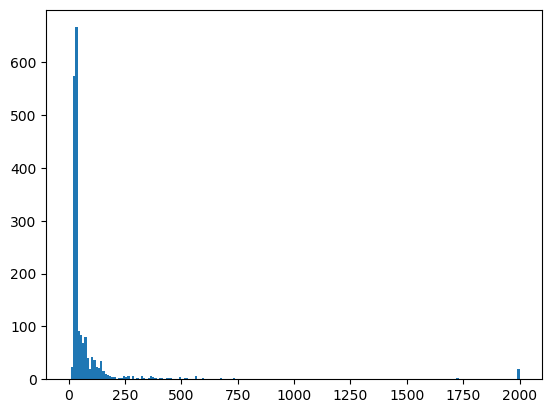}
    \caption{Raw spectral distribution,}
    \label{fig10:0}
  \end{subfigure}
  \begin{subfigure}[b]{0.45\textwidth}
    \includegraphics[width=\textwidth]{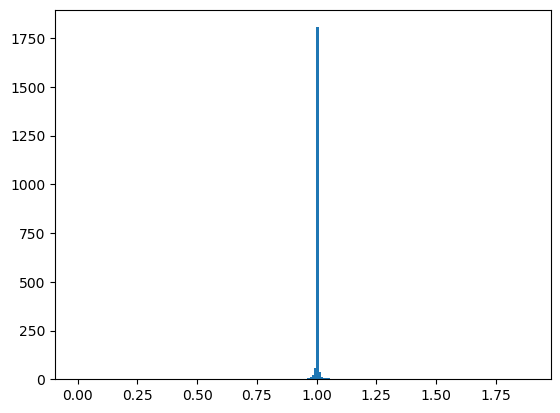}
    \caption{Normalized spectral distribution}
    \label{fig10:1}
  \end{subfigure}
  \caption{Positive rank 1 model , $p=0.05$}\label{pos1spec05}
\end{figure}
\begin{figure}
  \begin{subfigure}[b]{0.45\textwidth}
    \includegraphics[width=\textwidth]{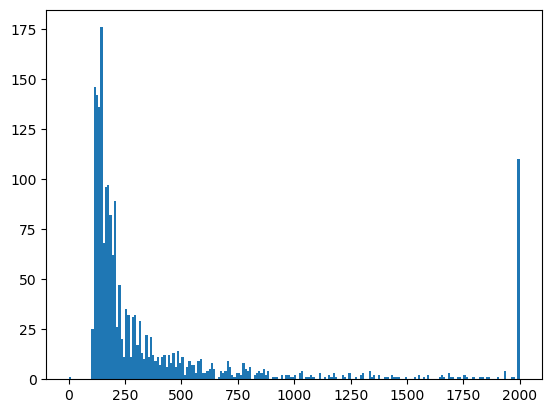}
    \caption{Raw spectral distribution,}
    \label{fig11:0}
  \end{subfigure}
  \begin{subfigure}[b]{0.45\textwidth}
    \includegraphics[width=\textwidth]{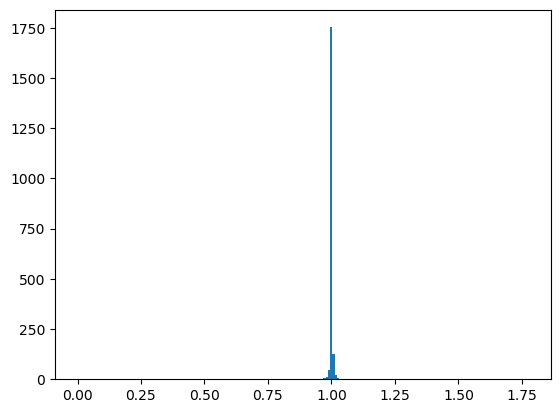}
    \caption{Normalized spectral distribution}
    \label{fig11:1}
  \end{subfigure}
  \caption{Positive rank 1 model , $p=0.2$}\label{pos1spec2}
\end{figure}
The raw distribution is interesting (there is a large spike at $n$), but what is \emph{more} interesting is that the normalized Laplacian has extreme concentration of eigenvalues at $1,$ completely unlike the Erd\"os-R\'enyi model. The standard deviation of the spectral distribution is (not surprisingly) much smaller, and it is also much less regular, the peak is also achieved for a far larger $p.$
\begin{figure}
  \begin{subfigure}[b]{0.45\textwidth}
    \includegraphics[width=\textwidth]{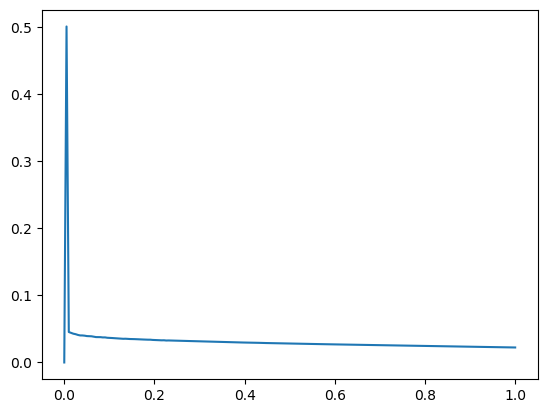}
    \caption{Entire graph}
    \label{fig12:0}
  \end{subfigure}
  \begin{subfigure}[b]{0.45\textwidth}
    \includegraphics[width=\textwidth]{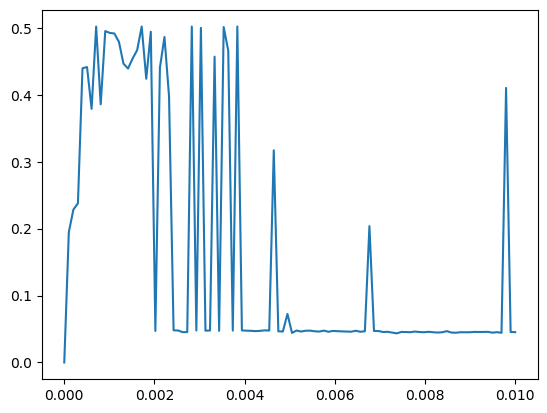}
    \caption{Zoom around maximum}
    \label{fig12:1}
  \end{subfigure}
  \caption{Standard deviation of spectral density for positive rank 1 model}\label{pos1std}
\end{figure}
\subsection{Wishart rank 1}
\begin{figure}
  \begin{subfigure}[b]{0.3\textwidth}
    \includegraphics[width=\textwidth]{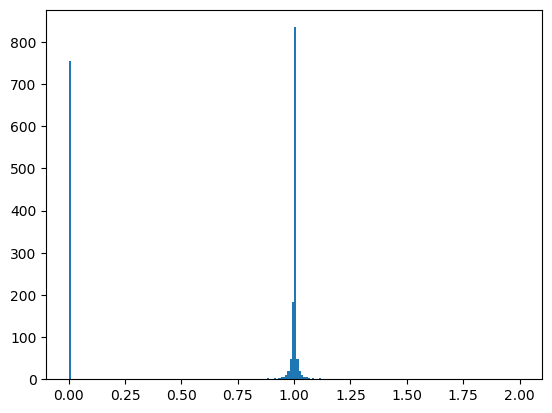}
    \caption{$p=0.05$}
    \label{fig13:0}
  \end{subfigure}
  \begin{subfigure}[b]{0.3\textwidth}
    \includegraphics[width=\textwidth]{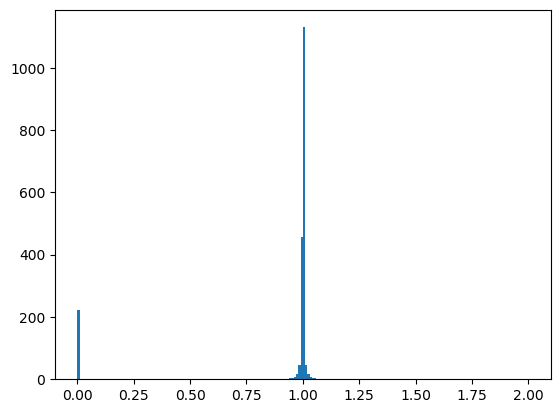}
    \caption{$p=0.2$}
    \label{fig13:1}
  \end{subfigure}
  \begin{subfigure}[b]{0.3\textwidth}
    \includegraphics[width=\textwidth]{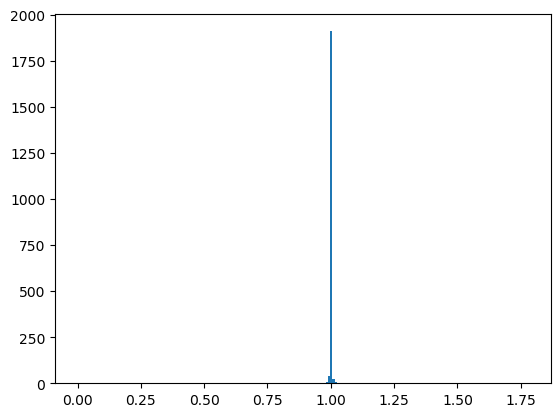}
    \caption{$p=0.6$}
    \label{fig13:2}
  \end{subfigure}
  \caption{Spectral density of Wishart rank 1 model}\label{wish1spec}
\end{figure}

\begin{figure}
  \begin{subfigure}[b]{0.45\textwidth}
    \includegraphics[width=\textwidth]{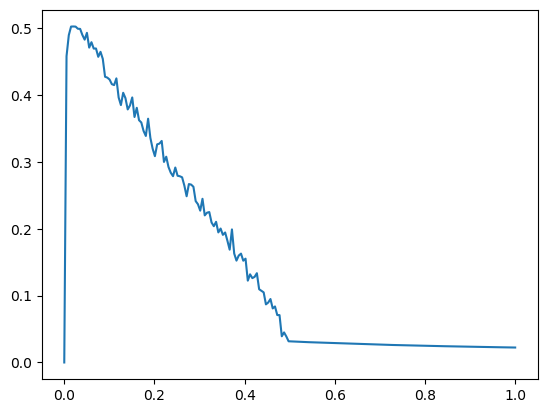}
    \caption{Entire graph}
    \label{fig14:0}
  \end{subfigure}
  \begin{subfigure}[b]{0.45\textwidth}
    \includegraphics[width=\textwidth]{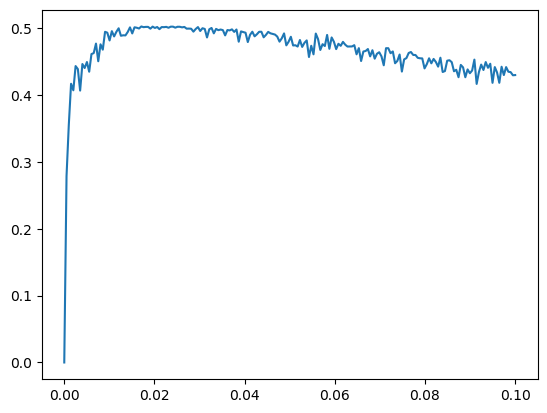}
    \caption{Zoom around maximum}
    \label{fig14:1}
  \end{subfigure}
  \caption{Standard deviation of spectral density for Wishart rank 1 model}\label{wish1std}
\end{figure}

The Wishart rank one graphs show essentially the same behavior as the positive rank one case, with a very tight concentration around $1,$ and rapid decay, but also a massive concentration at $0$ (indicating many connected components) for $p<0.5$ See Figure \ref{wish1spec}.
The standard deviation is quite different from the positive case - see Figure \ref{wish1std} - showing the usual phase transition at $p=0.5$
\subsection{Point Clouds} 
\begin{figure}
  \begin{subfigure}[b]{0.45\textwidth}
    \includegraphics[width=\textwidth]{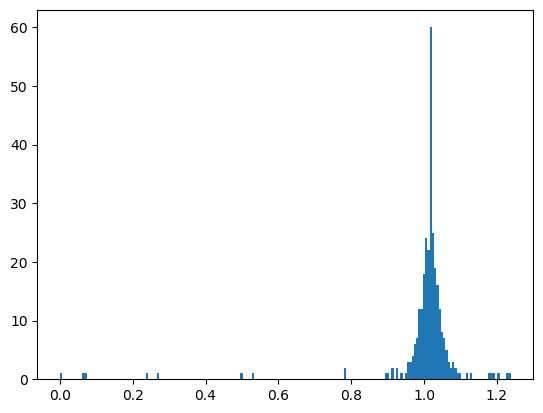}
    \caption{Noisy circle}
    \label{fig15:0}
  \end{subfigure}
  \begin{subfigure}[b]{0.45\textwidth}
    \includegraphics[width=\textwidth]{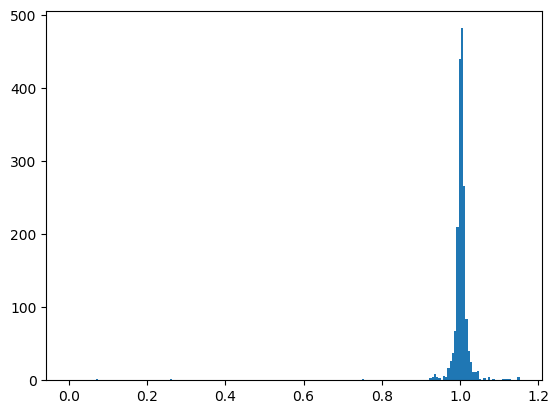}
    \caption{noisy torus}
    \label{fig15:1}
  \end{subfigure}
  \caption{Spectral density at $p=0.2$l}\label{noisyspec}
\end{figure}

\begin{figure}
  \begin{subfigure}[b]{0.45\textwidth}
    \includegraphics[width=\textwidth]{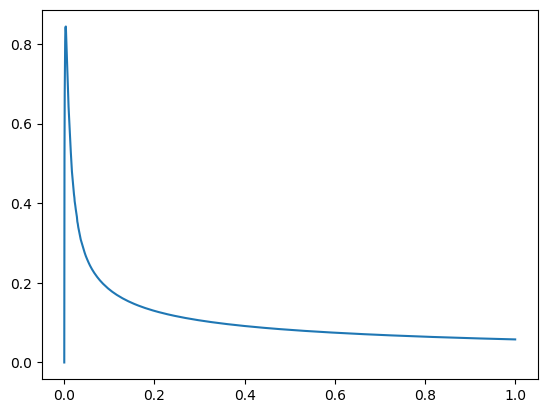}
    \caption{Entire graph}
    \label{fig16:0}
  \end{subfigure}
  \begin{subfigure}[b]{0.45\textwidth}
    \includegraphics[width=\textwidth]{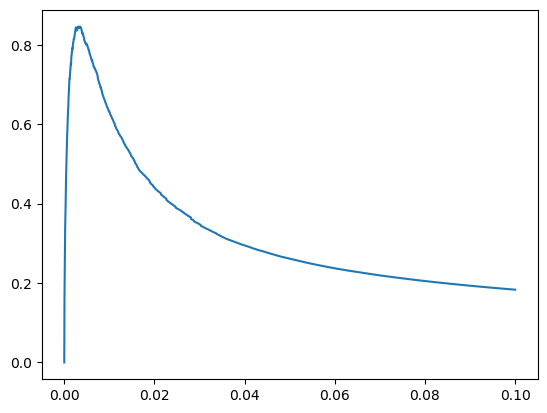}
    \caption{Zoom around maximum}
    \label{fig16:1}
  \end{subfigure}
  \caption{Standard deviation of spectral density for noisy circle}\label{circstd}
\end{figure}

\begin{figure}
  \begin{subfigure}[b]{0.45\textwidth}
    \includegraphics[width=\textwidth]{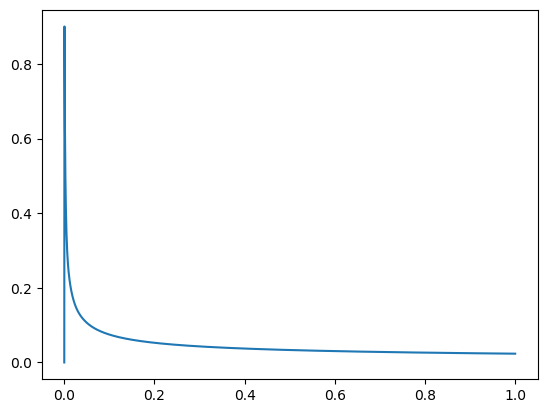}
    \caption{Entire graph}
    \label{fig17:0}
  \end{subfigure}
  \begin{subfigure}[b]{0.45\textwidth}
    \includegraphics[width=\textwidth]{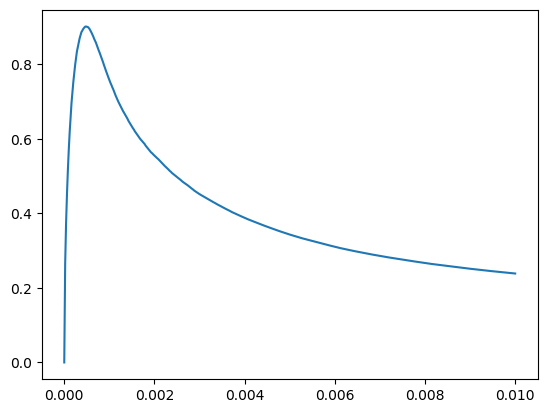}
    \caption{Zoom around maximum}
    \label{fig17:1}
  \end{subfigure}
  \caption{Standard deviation of spectral density for noisy torus}\label{torusstd}
\end{figure}

Here we look at the spectral distribution of the point clouds (noisy circle and noisy torus). It is evident that these are very close to the positive rank one matrices - the reader can judge for his or her own self. The bulk density at $p=0.2$ is in Figure \ref{noisyspec}
The evolution of standard deviation for the circle is given in Figure \ref{circstd}, for the torus in Figure \ref{torusstd}.
\bibliographystyle{alpha}
\bibliography{sample}

\end{document}